\documentclass[a4paper,12pt]{amsart}
\usepackage{latexsym}
\usepackage{amsfonts}
\usepackage{amsmath}
\usepackage{mathrsfs}
\usepackage{amsthm}
\usepackage{amssymb}

\newcommand\Proof{{\it Proof\/}}

\newtheorem{theorem}{Theorem}

\newcommand{\be}{\begin{eqnarray}}
\newcommand{\ee}{\end{eqnarray}}
\newcommand{\bestar}{\begin{eqnarray*}}
\newcommand{\eestar}{\end{eqnarray*}}
\newcommand{\ignore}[1]{}{}
\allowdisplaybreaks
\begin{document}

\vspace{1in}

\title[NON-ASYMPTOTIC RESULTS FOR CORNISH-FISHER EXPANSIONS]{\bf NON-ASYMPTOTIC RESULTS FOR CORNISH-FISHER EXPANSIONS}


\author[V.V. Ulyanov]{V.V. Ulyanov}
\address{V.V. Ulyanov\\
	Faculty of Computational Mathematics and Cybernetics\\
	Moscow State University \\
	Moscow, 119991, Russia\\
	and National Research University Higher School of Economics (HSE),
	Moscow, 101000, Russia
}
\email{vulyanov@cs.msu.ru}

\author[M. Aoshima]{M. Aoshima}
\address{M. Aoshima\\
	Institute of Mathematics\\
	University of Tsukuba\\ 
	Tsukuba, Ibaraki 305-8571, Japan
}
\email{aoshima@math.tsukuba.ac.jp}

\author[Y. Fujikoshi]{Y. Fujikoshi}
\address{Y. Fujikoshi\\
		Department of Mathematics\\
	Hiroshima University\\ 
		Higashi-Hiroshima, 739-8526, Japan
}
\email{fujikoshi\_y@yahoo.co.jp}

\thanks{This work was supported by RSCF grant No. 14–11–00364}


\begin{abstract}
We get the computable error bounds for generalized Cornish-Fisher expansions for quantiles of statistics 
provided that the computable error bounds for Edgeworth-Chebyshev type expansions for distributions of these statistics are known. The results are illustrated by examples.
\end{abstract}

\maketitle

\renewcommand{\refname}{References}

\section{Introduction and main results} 

 In statistical inference it is of fundamental importance to obtain the sampling distribution of statistics. However, we often encounter situations, where the exact distribution cannot be obtained in closed form, or even if it is obtained, it might be of little use because of its complexity. One practical way of getting around the problem is to provide reasonable approximations of the distribution function and its quantiles, along with extra information on their possible errors. It can be made with help of Edgeworth--Chebyshev and Cornish--Fisher type expansions. Recently the interest for Cornish--Fisher type expansions stirred up because of intensive study of VaR (Value at Risk) models in financial mathematics and financial risk management (see, e.g. \cite{Jash} and \cite{UVV:CFE}).
 
 Mainly, it is studied the asymptotic behavior of the expansions  mentioned above. It means that 
 accuracy of approximation for distribution of statistics or its quantiles is given as $O(\cdot)$, that is in the form of order with respect to some parameter(s) (usually, w.r.t. $n$ as a number of observations and/or $p$ as dimension of observations). In this paper we construct non-asymptotic error bounds, in other words -- computable error bounds, for Cornish--Fisher type expansions, that is for an error of approximation  we prove upper bounds with closed-form dependence on $n$ and/or $p$ and, perhaps, on some moment characteristics of observations. We get these bounds under condition that similar non-asymptotic  results are already known for accuracy of approximation of distributions of statistics by Edgeworth--Chebyshev type expansions. 
 
 

 Let $X$ be a univariate random variable with a continuous distribution function $F$. For $\alpha :\, 0<\alpha <1$, there exists $x$ such that $F(x)=\alpha$, which is called the (lower) 100$\alpha \%$
 \textit{point} of $F$. If $F$ is strictly increasing, the inverse function $F^{-1}(\cdot)$  is well defined and the 100$\alpha \%$
 point is uniquely determined.
 We also speak of ``quantiles'' without reference to particular values of $\alpha$ meaning the values given by    $F^{-1}(\cdot)$.
 %
 Even in the general case, when $F(x)$ is not necessarily continuous nor is it strictly increasing, we can define its inverse function by formula
 $$
 F^{-1}(u) =\inf\{x; F(x) > u\}.
 $$
 This is a right-continuous nondecreasing function defined on the interval $(0, 1)$ and $F(x_0) \ge u_0$ if $x_0 =F^{-1}(u_0)$.
 \par
 Let $F_n(x)$ be a sequence of distribution functions and  let
 each $F_n$ admit  the Edgeworth-Chebyshev type  expansion (ECE) in the powers of $\epsilon=n^{-1/2}$ or $n^{-1}$:
 \begin{equation}
 \begin{split}
 F_n(x) &= G_{k,n}(x)+R_k(x) ~~~ {\rm with}\,\, R_k(x) = O(\epsilon^k) ~~~ {\rm and}
 \\
 G_{k,n}(x) &= G(x)+\bigl\{
 \epsilon a_1(x) + \dots + \epsilon^{k-1}a_{k-1}(x)
 \bigr\}
 g(x),
 \end{split}
 \label{eq5-6-1}
 \end{equation}
 where $g(x)$ is a density function of the limiting distribution function $G(x)$.
 An important approach to the problem of approximating the quantiles of $F_n$ is to use their asymptotic relation to those of $G$'s.
 Let $x$ and $u$ be the corresponding quantiles of $F_n$ and $G$, respectively. Then we have
 %
 \begin{equation}
 F_n(x)=G(u).
 \label{eq5-6-3}
 \end{equation}
 Write $x(u)$ and $u(x)$ to denote the solutions of \eqref{eq5-6-3} for $x$ in terms of $u$ and $u$ in terms of $x$, respectively [i.e. $u(x)=G^{-1}(F_n(x))$ and $x(u)=F_n^{-1}(G(u))$].
 Then we can use the ECE \eqref{eq5-6-1}   to obtain formal solutions $x(u)$ and $u(x)$ in the form
 %
 \begin{equation}
 x(u)=u+\epsilon b_1(u)+\epsilon^2 b_2(u)+ \cdots ~
 \label{eq5-6-4}
 \end{equation}
 \vskip -8pt
 \noindent
 and
 \vskip -8pt
 %
 \begin{equation}
 u(x)=x+\epsilon c_1(x)+\epsilon^2 c_2(x)+ \cdots.
 \label{eq5-6-5}
 \end{equation}
 \par
 Cornish and Fisher in \cite{CF} and \cite{FC}  obtained the first few terms of these expansions when $G$ is the standard normal distribution function (i.e., $G= \Phi$).
 Both \eqref{eq5-6-4} and \eqref{eq5-6-5} are called the {\it Cornish--Fisher expansions, (CFE)}. Concerning CFE for random variables obeying limit laws from the family of Pearson distributions see, e.g., \cite{Bo}. 
 Hill and Davis in \cite{HD} gave a general algorithm for obtaining each term of CFE when $G$ is an analytic function.
 
 Usually the CFE are applied in the following form with $k= 1,\,2$ or $3$:
 \begin{equation}
 x_{k}(u)=u+ \sum_{j=1}^{k-1}\epsilon^j b_j(u)+ \hat{R}_k(u) ~~~ {\rm with}\,\, \hat{R}_k(u) = O(\epsilon^k).
 \label{12Fn36}
 \end{equation}
It is known (see, e.g., \cite{UVV:CFE}) how to find the explicit expressions for $b_1 (u)$ and $b_2 (u)$ as soon as we have \eqref{eq5-6-1}. By Taylor's expansions for $G$, $g$, and $a_1$, we obtain
\begin{equation}
\begin{split}
b_1 &=
-a_1(u),
 \\
b_2 &=
\frac{1}{2} \{ g'(u)/g(u) \} a_1^2(u)
-a_2(u) + a_1'(u) a_1(u), 
\end{split}
\label{eq5-6-8}
\end{equation}
provided that $g$ and $a_1$ are smooth enough functions.

In the following Theorems we show how $x_{k}(u)$ from \eqref{12Fn36} could be expressed in terms of $u$. Moreover, we show what kind of bounds we can get for $\hat{R}_k(x)$ as soon as we have some bounds for ${R}_k(x)$ from \eqref{eq5-6-1}.
 
\begin{theorem}
	\label{th5-6-2}
	Suppose that for the distribution function of a statistic $U$ we have
	\begin{equation}
	F(x) \equiv \Pr\{U \le x\}=G(x)+R_1(x), 
	\label{eq5-6-20}
	\end{equation}
	where for remainder term
	$R_1(x)$
	there exists a constant $c_1$ such that
	$$
	|R_1(x)| \le d_1 \equiv c_1 \epsilon.\,\,
	$$
	Let $x_{\alpha}$ and $u_{\alpha}$ be the upper $100 \alpha \%
	$
	points of $F$ and $G$ respectively, that is
	\begin{equation}
	\Pr\{U \le x_{\alpha}\}=G(u_{\alpha})=1-\alpha. 
	\label{eq5-6-21}
	\end{equation}
	Then for any $\alpha$ such that $1 - c_1 \epsilon > \alpha > c_1 \epsilon > 0$ we have
	\begin{description} \setlength{\itemsep}{3pt}
		\item[\hspace{2pt}{\rm (i)}] \hspace{5pt}$ {u_{\alpha+d_1} \le x_{\alpha} \le  u_{\alpha-d_1}}$.
		\item[{\rm (ii)}] \hspace{3pt}
		${|x_{\alpha}-u_{\alpha}| \le c_1 \epsilon /g(u_{(1)})}$, where $g$ is the density function of the limiting distribution
		$G$ and		
	\end{description}	
	$$
	g(u_{(1)})= \min_{u \in [u_{\alpha+d_1}, u_{\alpha-d_1}]} g(u). \hspace{64pt}~
	$$
\end{theorem}

\begin{theorem}
	\label{th5-6-2}
	In the notation of Theorem 1 we assume that
	\begin{equation}
	F(x) \equiv \Pr\{U \le x\}=G(x)+\epsilon g(x)a(x) + R_2(x),\nonumber
	\label{eq5-6-2}
	\end{equation}
	where for remainder term
	$R_2(x)$
	there exists a constant $c_2$ such that
	$$
	|R_2(x)| \le d_2 \equiv c_2 \epsilon^2.
	$$
	Let $T=T(u)$ be a monotone increasing transform such that
	$$\Pr\{T(U) \le x\} = G(x) + \tilde R_2(x) \,\,\,\,\text{with} \,\,\,\,|\tilde R_2(x)| \leq \tilde d_2 \equiv \tilde c_2 \epsilon^2.
	$$
	Let $\tilde x_{\alpha}$ and $u_{\alpha}$ be the upper $100 \alpha \%
	$
	points of $\Pr\{T(U) \le x\}$ and $G$, respectively.
	Then for any $\alpha$ such that 
	$$
	1 - \tilde c_2 \epsilon^2 > 
	\alpha > \tilde c_2 \epsilon^2 > 0,
	$$ 
	we have
	\begin{equation}
 	{|\tilde x_{\alpha}-u_{\alpha}| 
 		\le \tilde c_2 \epsilon^2 /g(u_{(2)})},
 	\label{Th2.1}
 	 \end{equation}
 	 where
	%
	$$
	g(u_{(2)})= \min_{u \in [u_{\alpha+ \tilde d_2}, u_{\alpha- \tilde d_2}]} g(u). \hspace{64pt}~
	$$
\end{theorem}

\begin{theorem}
	\label{th5-6-2} We use the notation of Theorem 2.
	Let $b(x)$ be a function inverse to $T$, i.e. $b(T(x))=x$. Then $	x_\alpha = b(\tilde x_\alpha)$
	and for $\alpha$ such that $1 - \tilde c_2 \epsilon^2 > \alpha > \tilde c_2 \epsilon^2$ we have
	\begin{equation}
	|x_\alpha - b(u_\alpha)| \leq \tilde c_2 \,  \frac{|b'(u^*)|}{g(u_{(2)})}\,\epsilon^2,
	\label{Th3.1}
	\end{equation}
	where $$|b'(u^*)| = \max_{u\in[u_{\alpha+\tilde d_2},u_{\alpha-\tilde d_2}]} |b'(u)|.$$
	Moreover, 
	\begin{equation}
	b(x) = x - \epsilon a(x) + O(\epsilon^2).
	\label{Th3.0}
	\end{equation}
\end{theorem}

{\bf Remark 1.} The main assumption of the Theorems is that for distributions of statistics  and for distributions of transformed statistics we have some approximations with computable error bounds. There are not many papers with this kind of non-asymptotic results because it requires technique which is different from the asymptotic results methods (cf., e.g., \cite{FUS2005} and \cite{WFU2014}). In series of papers \cite{FU2006}, \cite{FU2006b}, \cite{FUS2005}, \cite{FUS2005b}, \cite{GUF2013}, \cite{UCF2006}, \cite{UFS1999}, \cite{UWF2006}  we got non-asymptotic results for wide class of statistics including multivariate scale mixtures and MANOVA tests. We considered as well the case of high dimensions, that is the case when the dimension of observations and sample size are comparable. The results were included in the book \cite{UVV:FUS}. See also \cite{EA2}.

{\bf Remark 2.} The results of Theorems 1--3 could not be extended to the whole range  of $\alpha\in (0,1)$. It follows from the fact that the Cornish-Fisher expansion does not converge uniformly in $0<\alpha<1.$ See corresponding example in 
Section 2.5 of \cite{Hall}.

{\bf Remark 3.} In Theorem 2 we required the existence of a monotone increasing transform $T(z)$ such that distribution of transformed statistic $T(U)$ is approximated by some limit distribution $G(x)$ in better way than the distribution of original statistic $U$. We call this transformation $T(z)$ the   {\it Bartlett type correction}. See corresponding examples in Section 3. 

{\bf Remark 4.} According to \eqref{Th3.1} and 
 \eqref{Th3.0} the function $b(u_{\alpha})$ in Theorem 3 could be considered as an  ''asymptotic expansion'' for $x_{\alpha}$ up to order $O(\epsilon^2).$

\section{Proofs of main results}

\Proof {\it of Theorem 1.} By the mean value theorem,
$$
|G(x_{\alpha})-G(u_{\alpha})| \ge
|x_{\alpha}-u_{\alpha}| \min_{0<\theta<1}
g(u_{\alpha}+\theta(x_{\alpha}-u_{\alpha})).
$$
From \eqref{eq5-6-20} and the definition of $x_{\alpha}$ and $u_{\alpha}$ in \eqref{eq5-6-21},
we get
%
\begin{align*}
|G(x_{\alpha})-G(u_{\alpha})|
&= \left|G(x_{\alpha}) - \Pr\{U \le x_{\alpha}\}\right| \\
&= |R_1(x_{\alpha})| \le d_1.
\end{align*}
Therefore,
%
\begin{equation}
|x_{\alpha}-u_{\alpha}|
\le \frac{d_1}{\min_{0<\theta<1}
	g(u_{\alpha}+\theta(x_{\alpha}-u_{\alpha})}.
\label{eq5-6-22}
\end{equation}
On the other hand, it follows from \eqref{eq5-6-20} that
%
\begin{align*}
G(x_{\alpha})&= 1-\alpha-R_1(\alpha) \\
&\le 1-(\alpha -d_1)=G(u_{\alpha -d_1}).
\end{align*}
This implies that $x_{\alpha} \le u_{\alpha -d_1}$.
Similarly, we have $u_{\alpha +d_1} \le x_{\alpha}$.
Therefore, we proved Theorem 1 (i).

It follows from Theorem 1 (i) that  
%
$$
\min_{u \in [u_{\alpha+d_1}, u_{\alpha-d_1}]} g(u)
\le \min_{0<\theta<1}
g(u_{\alpha}+\theta(x_{\alpha}-u_{\alpha})).
$$
Thus, using \eqref{eq5-6-22} we get statement of Theorem 1 (ii).
\\
\\
\Proof {\it of Theorem 2.} It is easy to see that it is sufficient to apply Theorem 1 (ii) to the transformed statistic $T(U).$
\\
\\
\Proof {\it of Theorem 3.} Using now \eqref{Th2.1} and the mean value theorem we obtain
\begin{equation} 
\tilde x_{\alpha}-u_{\alpha} = b^{-1}(x_\alpha) - b^{-1}(b(u_\alpha)) = (b^{-1})'(x^*) \big(x_\alpha - b(u_\alpha)\big)\,,
\label{Th3.2}
\end{equation}
where $x^*$ is a point on the interval $\big(\min\{x_\alpha , b(u_\alpha)\} \,,\, \max\{x_\alpha , b(u_\alpha)\}\big)$. \\
By Theorem 1 (i) we have
$$
{u_{\alpha+\tilde d_2} \le \tilde{x}_{\alpha} \le  u_{\alpha-\tilde d_2}}.
$$
Therefore, for $	x_\alpha = b(\tilde x_\alpha)$ we get
\begin{equation}
\big(\min\{b^{-1}(x_\alpha), u_{\alpha}\} \,,\, \max\{b^{-1}(x_\alpha), u_{\alpha}\}\big) \subseteq \big(u_{\alpha+\tilde d_2} \,,\, u_{\alpha-\tilde d_2} \big).
\label{Th3.3}
\end{equation}
Since by properties of derivatives of inverse functions
 $$
 (b^{-1})'(z) = 1/\, b'(b^{-1}(z))= 1/b'(y)
 $$ 
 for $z=b(y)$, the relations \eqref{Th3.2} and \eqref{Th3.3} imply \eqref{Th3.1}. 
 
 Representation \eqref{Th3.0} for $b(x)$ follows from \eqref{eq5-6-8} and \eqref{Th3.1}.
 
\section{Examples}
 
In \cite{UF2001} we gave sufficient conditions for transformation $T(x)$ to be the  Bartlett type correction (see Remark 3 above) for wide class of statistics $U$ allowing the following represantion 
\begin{equation}
\Pr\{U \le x\}=G_q(x)+\frac{1}{n}\,\sum_{j=0}^{k}a_j\,G_{q+2j}(x)  +  R_{2k},
\label{Ex1}
\end{equation}
where $R_{2k} = O(n^{-2})$ and  $G_q(x)$ is the distribution function of chi-squared distribution with $q$ degrees of freedom and coefficients $a_j$'s satisfy the relation $\sum_{j=0}^{k}a_j = 0$. Some
examples of the statistic $U$ are as follows: for $k=1,$ the likelihood ratio test statistic; for $k=2,$ the Lawley-Hotelling trace criterion and the Bartlett-Nanda-Pillai trace 
criterion, which are test statistics for multivariate linear hypothesis under normality; for $k=3,$ the score test statistic  and  Hotelling's $T^{2}$-statistic under nonnormality. The results of \cite{UF2001} were extended in \cite{EA1} and \cite{EA2}. 

In \cite{FUS2005} we were interested in the null distribution of Hotelling's generalized $T^2_0$
 statistic defined by
$$
T^2_0 = n\,\mbox{tr} S_h\,S^{-1}_e,
$$
where $S_h$ and $S_e$ are independently distributed as Wishart distributions $W_p(q, I_p)$ and
$W_p(n, I_p)$ with identity operator $I_p$ in $\mathbb{R}^p$, respectively. In Theorem 4.1 (ii) in \cite{FUS2005} we proved \eqref{Ex1} for all $n\geq p$ 
with $k=3$ and  computable error bound:  
\begin{eqnarray*}
	&&|{\Pr}(T_0^2 \leq x)-G_r(x)-\frac{r}{4n}\{(q-p-1)G_r(x)  \\
	&& \quad -2qG_{r+2}(x)+(q+p+1)G_{r+4}(x)\}|  \\
	&& \ \quad \leq \frac{c_{p,q}}{n^2}, 
\end{eqnarray*}
where $r=pq$ and for constant $c_{p,q}$ we gave expicit formula with dependence on $p$ and $q$.

 Therefore, according to \cite{UF2001} we can take in this case the Bartlett type correction $T(z)$ as
 $$
 T(z)=\frac{a-1}{2b} +
 \sqrt{\left(\frac{a-1}{2b}\right)^2+\frac{z}{b}},
 $$
 where
 $$
 a=\frac{1}{2n}p(q-p-1),$$
 $$
 b=\frac{1}{2n}p(q+p+1)(q+2)~{-1}.
 $$
 It is clear that $T(z)$ is invertable and we can apply Theorem 3.
 
 Other examples and numerical calculations and comparisons of approximation accuracy  see in \cite{EA1} and \cite{EA2}.
 
 One more example is connected with sample correlation coefficient. Let $\vec{X}=(X_1,...,X_n)^T$, and $\vec{Y}=(Y_1,...,Y_n)^T$ be two vectors from an $n$-dimensional normal distribution $N(0, I_n)$ with zero mean, identity covariance matrix $I_n$ and
 the  sample correlation coefficient
 \begin{equation}\label{g1c}
 R=R(\vec{X},\vec{Y})=\frac{\sum\nolimits_{k=1}^n\,X_k\,Y_k}{\sqrt{\sum\nolimits_{k=1}^n\,X_k^2 \,\,\, \sum\nolimits_{k=1}^n\,Y_k^2}}\,.\nonumber 
 \end{equation}
 In \cite{GUF2013} it was proved for $n \geq 7$ and $N=n-2.5$:  
 \begin{equation}\label{t3}
 \sup\nolimits_x \left| \Pr\Big(\sqrt{N}\,R \leq x\Big) - \Phi(x) - \frac{x^3\,\varphi(x) }{4\,N}  \right|\leq \frac{B_n}{N^2}, \nonumber
 \end{equation}
 with $B_n \leq 2.2.$
 It is easy to see that we can take $T(z)$ as the Bartlett type correction in the form
   $T(z) =  z+z^3/(4N)$.  Then 
   the inverse function $b(z)=T^{-1}(z)$ is defined by formula 
 \begin{eqnarray}\label{g72a}
 b(z) &=& \Big(2\,N\,z+\sqrt{(2Nz)^2 + (4N/3)^3}\Big)^{1/3}\nonumber\\
 &&\,\,-\Big(-\,2\,N\,z+\sqrt{(2Nz)^2 + (4N/3)^3}\Big)^{1/3} \nonumber\\
 &=& z - \frac{z^3}{4\,N} + \frac{3\,z^5}{16\,N^2} + O(N^{-3}). \nonumber 
 \end{eqnarray}
 Now we can apply Theorem 3.

\end{document}